\documentclass[11pt, leqno]{article}
\usepackage{amsmath,amsfonts,amssymb,wasysym,mathrsfs}
\usepackage[latin1]{inputenc}
\usepackage{lmodern}
\usepackage{shapepar}
\usepackage{graphicx}
\usepackage[T1]{fontenc}
\usepackage{frcursive}
\usepackage{hyperref}
\usepackage{ifpdf} 

\usepackage{color}
\DeclareMathAlphabet{\mathpzc}{OT1}{pzc}{m}{it}
\newcommand{\R}{{\mathbb R}}
\newenvironment{frcseries}{\fontfamily{frc}\selectfont}{}

\newcommand{\be}[1]{\begin{equation}\label{#1}}
\newcommand{\ee}{\end{equation}}

\newcommand{\mathfrc}[1]{\text{\textfrc{#1}}}
\newcommand{\textfrc}[1]{{\frcseries#1}}

\newcommand{\prf}{\par\smallskip\noindent{\sl Proof. \/}}
\newcommand{\finprf}{\unskip\null\hfill$\;\square$\vskip 0.3cm}

\newtheorem{theorem}{Theorem}[section]
\newtheorem{lemma}{Lemma}[section]
\newtheorem{corollary}[theorem]{Corollary}

\newtheorem{remark}[theorem]{Remark}
\newtheorem{definition}{Definition}[section]

\numberwithin{equation}{section}
\newcommand{\nc}{\normalcolor}

\def\qed{\,\unskip\kern 6pt \penalty 500
\raise -2pt\hbox{\vrule \vbox to8pt{\hrule width 6pt
\vfill\hrule}\vrule}\par}
\definecolor{darkblue}{rgb}{0.05, .05, .65}
\definecolor{darkgreen}{rgb}{0.1, .65, .1}
\definecolor{darkred}{rgb}{0.8,0,0}
\begin{document}
\title{\textbf{Symmetrization for fractional nonlinear elliptic problems}\\[7mm]}

\author{\Large Vincenzo Ferone\footnote{Dipartimento di Matematica e Applicazioni ``Renato Caccioppoli'', Universit\`a degli Studi di
Napoli Federico II, 80143  Italia. \    E-mail: {\tt vincenzo.ferone@unina.it}}  \ and \ Bruno Volzone \footnote{Dipartimento di Scienze e Tecnologie, Universit\`a degli Studi di
Napoli ``Parthenope'', 80143  Italia. \    E-mail: {\tt bruno.volzone@uniparthenope.it}}
\\[0.5cm]
(\emph{Dedicated to Juan Luis V\'azquez, a maestro and friend, for his 75th birthday})}

\date{} 

\maketitle

\begin{abstract}
In this note we prove a new symmetrization result, in the form of mass concentration comparison, for solutions of nonlocal nonlinear Dirichlet problems involving fractional $p$ Laplacians. Some regularity estimates of solutions will be established as a direct application of the main result.
\end{abstract}

\setcounter{page}{1}
\section{Introduction} \label{Sec1}
The main goal of this note is to provide new results concerning the application of symmetrization methods in the context of \emph{nonlocal, nonlinear} elliptic problems. In particular, we will focus in getting new estimates for solutions to the elliptic problem
\begin{equation} \label{mainpp}
\left\{
\begin{array}[c]{lll}%
\left( -\Delta_{p}\right)^{s}u=f & & \text{in }%
\Omega,\\
\\
u=0 & & \text{on }\R^{N}\setminus\Omega.
\end{array}\right.
\end{equation}
The operator $\left( -\Delta_{p}\right)^{s}$ is the so called \emph{fractional} \emph{p-Laplacian} and is defined for $s\in (0,1)$ and $p>1$ by means of the singular integral formula
\[
(-\Delta_p)^su(x)=\gamma(N,s,p)\;\text{P.V.} \int_{\R^N}\frac{|u(x)-u(y)|^{p-2} (u(x)-u(y))}{|x-y|^{N+sp}}dy,
\]
being $\gamma(N,s,p)$ a suitable normalization constant, whose value is specified as (see \cite[Lemma 5.1]{CDVAZ})
\begin{equation*}\label{constant}
\gamma(N,s,p)=\frac{sp\,2^{2s-2}(1-s)}{\pi^{\frac{N-1}{2}}}\frac{\Gamma\left(\frac{N+sp}{2}\right)}{\Gamma\left(\frac{p+1}{2}\right)\Gamma(2-s)}.
\end{equation*}
In this note we will consider only the \emph{degenerate} case $p>2$. More details regarding the correct functional spaces for the source term $f$ and the corresponding solution $u$ will be given in Subsection \ref{Functionspac}. Moreover, we assume $N\geq1$ while the ground set $\Omega\subset \R^{N}$ will be always assumed bounded with Lipschitz boundary.\\[0.2pt]

\noindent Inspired by the techniques of the recent work \cite{FeroneVolzone}, in this note we derive some symmetrization estimates for solutions $u$ to problem \eqref{mainpp}, in the form of \emph{mass concentration comparison}. For the sake of completeness, we recall the main result obtained in \cite{FeroneVolzone}, which was also established in \cite{VazVolSire} through the Caffarelli-Silvestre extension theorem. If we consider the linear, nonlocal Dirichlet problem
\begin{equation} \label{symmline}
\left\{
\begin{array}[c]{lll}%
\left( -\Delta\right)^{s}u=f & & \text{in }%
\Omega,\\
\\
u=0 & & \text{on }\R^{N}\setminus\Omega.
\end{array}\right.
\end{equation}
then the \emph{worst} radial problem, that is the problem whose solution $v$ is the larger one in the class of problems of the type \eqref{mainpline} with domains $\Omega$ of fixed measure and corresponding data $f$, solves the \emph{symmetrized} problem
\begin{equation} \label{mainpline}
\left\{
\begin{array}[c]{lll}%
\left( -\Delta\right)^{s}v=f^{\#} & & \text{in }%
\Omega^{\#},\\
\\
v=0 & & \text{on }\R^{N}\setminus\Omega^{\#},
\end{array}\right.
\end{equation}
being $\Omega^{\#}$ the ball centred at the origin with volume $|\Omega|$: indeed, what we have is that
\begin{equation}\label{comparis}
u^{\#}\prec v
\end{equation}
where $u^{\#}$ is the \emph{Schwarz decreasing rearrangement of $u$ }and $\prec$ is the order relation in the form of mass concentration comparison (see Section \ref{Sec2} for precise definitions) . The machinery used in \cite{FeroneVolzone} consists in choosing a suitable truncation function in \eqref{symmline} and a convenient Riesz rearrangement inequality, which allows to obtain an integral estimate involving the Schwarz rearrangement $u^\#$ of $u$ and the datum $f$. Therefore, a quite subtle and technical work is made to reinterpret the latter estimate as a comparison between some integral mean functions of $u^\#$ and $v$ on balls, which in turn implies \eqref{comparis} by means of maximum principle arguments. Finally, an entire Section in \cite{FeroneVolzone} is dedicated to show that estimate \eqref{comparis} is \emph{optimal} for $s\in (0,1)$, in the sense that no pointwise comparison is achieved, while Talenti's classical result is recovered in the limit for $s\rightarrow1$. We refer the interested to \cite{FeroneVolzone} for all the details and for an exhaustive list of references concerning symmetrization results.\\\\[0.2pt]
\noindent As it was mentioned above, we wish to apply the \emph{direct} methods of \cite{FeroneVolzone} (\emph{i.e.} which does not use the extension theorem of \cite{Caffarelli-Silvestre}) to the nonlinear context \eqref{mainpp}.  For the \emph{local} context $s\rightarrow1$ the main reference is with no doubt \cite{talNON}, but several extentions to the case of more general local nonlinear problems can be found (see, for example, \cite{BFM}, \cite{AFT}, \cite{fermess}, \cite{cianchimaz}). In particular, in \cite{talNON} it is shown that if $u$ solves the nonlinear problem
\begin{equation*} \label{mainp}
\left\{
\begin{array}[c]{lll}%
-\Delta_{p}u=f & & \text{in }%
\Omega,\\
\\
u=0 & & \text{on }\partial\Omega
\end{array}\right.
\end{equation*}
being $-\Delta_{p}$ the classical $p$-Laplacian, and $v$ is the solution to
\begin{equation*} \label{mainp}
\left\{
\begin{array}[c]{lll}%
-\Delta_{p}v=f^{\#} & & \text{in }%
\Omega^{\#},\\
\\
v=0 & & \text{on }\partial\Omega^{\#}
\end{array}\right.
\end{equation*}
then a pointwise estimate is possible, namely we have
\[
u^{\#}\leq v \quad\text{in }\Omega^{\#}.
\]
Such a result is essentially based on the techniques of \cite{Talenti1}, where a basic importance is assigned to the explicit form fo the solution $v$, obtained by solving a nonlinear radial ODE. In particular,  the following estimate is obtained (with $r=|x|$)
\[
|\nabla u^\#(r)|\leq \frac{1}{r^{\frac{N-1}{p-1}}}\left(\int_{B_{r}}f^{\#}\right)^{\frac{1}{p-1}}
\]
but an explicit computation gives that the right-hand side equals $|\nabla v|$ and the result follows. \\ There is no need to say that when $\Omega$ is a ball, such ODE approach is not appliable to problems of the type \eqref{mainpp}, thus the nonlocal nature of the problem affects again the style of the technical approach.  Therefore our main goal will be to compare the solution $u$ to  \eqref{mainpp} with the solution $v$ of a suitable nonlocal \emph{radial} problem posed on the ball $\Omega^{\#}$, see problem \eqref{eq.symmetric}. Such choice of the symmetrized problem is justified by the appearance of a nonlinear integral estimate for the solution $u$, to which a maximum principle argument (in the spirit of the one derived in \cite{FeroneVolzone})  does not seem to apply in a direct form. In any case, ours is the first symmetrization result in the literature for nonlocal nonlinear problem and will allow to find $L^{p}$ regularity estimates for $u$ as a quite direct consequence. Finally, this note focuses only on the case $p\ge2$, because we shall need such convexity requirement for the application of a convenient Riesz's general rearrangement inequality.
\\[0.2pt]

\noindent The paper is organized as follows. Section \ref{Sec2} is entirely devoted to various preliminaries which will be used throughout the text: symmetrization tools will be shortly introduced and the suitable functional context will be settled. In Section \eqref{Main} we state the main comparison Theorem \ref{maint} and the $L^{p}$ regularity estimates of Corollary \eqref{regularity}. In Section \ref{Proofs} we will prove all the stated results. In Section \ref{Open problems} we illustrate open problems and numerical studies.

\section{Preliminaries and notation} \label{Sec2}
For the sake of completeness, we shortly collect here some preliminary results regarding symmetrization, the functional spaces related to the main problem \ref{mainpp} and hypergeometric functions.

\subsection{Rearrangements and symmetrization}\label{RearSym}

We recall here the basic definitions concerning Schwarz symmetrization and some related fundamental properties. Readers who wish to find more more details
of the theory are addressed to the classical monographs \cite{Hardy}, \cite{Bennett},  \cite{Kesavan}, \cite{Bandle} or to the paper
\cite{Talentirearrinv}.

A measurable real function $f$ defined on $\R^{N}$ is called \emph{radially symmetric} (or \emph{radial}) if there is a function
$\widetilde{f}:[0,\infty)\rightarrow \R$ such that $f(x)=\widetilde{f}(|x|)$ for all $x\in \R^{N}$. We will often write $f(x)=f(r)$,
$r=|x|\ge0$ for such functions by abuse of notation. We say that $f$ is \emph{rearranged} if it is radial, nonnegative and $\widetilde{f}$ is a
right-continuous, non-increasing function of $r>0$. A similar definition can be applied for real functions defined on a ball
$B_{R}(0)=\left\{x\in\R^{N}:|x|<R\right\}$.


Let $f$ be a real measurable function on $\R^N$. If $f$ is such that
its
\emph{distribution function} $\mu_{f}$ satisfies%
\begin{equation}\label{distribution}
\mu_{f}(  t)  :=\left\vert \left\{  x\in\\R^{N}:\left\vert f\left(
x\right)  \right\vert >t\right\} \right\vert<+\infty, \qquad\text{for every }t
>0,
\end{equation}
we define the \emph{one dimensional decreasing rearrangement} of $f$ as%
\[
f^{\ast}\left(  \sigma\right)  =\sup\left\{ t\geq0:\mu_{f}\left(  t\right)
>\sigma\right\}  \text{ , }\sigma>0.
\]
If $f$ is a real measurable function on an open set $\Omega\subset\R^N$ we extend $f$  as the  zero function in $\R^N\backslash\Omega$ and we define the one dimensional decreasing rearrangement of $f$ as the rearrangement of such an extension. This means that $f^{\ast}(\sigma)=0$ for $\sigma\in[|\Omega|,\infty)$. From the above definition it follows
that $\mu_{f^{\ast}}=\mu_{f}$ (\emph{i.e.,} $f$ and $f^{\ast}$ are equi-distributed) and $f^{\ast}$ is exactly  the \emph{generalized right
inverse function} of
$\mu_{f}$.
Furthermore, if $\Omega^{\#}$ is the ball of $%
\mathbb{R}
^{N}$ centered at the origin having the same Lebesgue measure as $\Omega$ ($\Omega^{\#}=\R^N$ if $|\Omega|=+\infty$), denoting by
\begin{equation*}
\omega_N=\frac{\pi^{N/2}}
{\Gamma\left(\frac N2+1\right)}
\end{equation*}
the measure of the unit ball in $\R^N$,
we define the
function
\[
f^{\#}\left(  x\right)  =f^{\ast}(\omega_{N}\left\vert x\right\vert
^{N})\text{ \ , }x\in\Omega^{\#},
\]
that will be called \emph{radially decreasing rearrangement}, or \emph{Schwarz decreasing rearrangement}, of $f$. We easily infer that $f$ is rearranged if and
only if
$f=f^{\#}$.

As a simple  consequence of the definition is that rearrangements preserve
$L^{p}$
norms, that is, for all $p\in[1,\infty]$
\[
\|f\|_{L^{p}(\Omega)}=\|f^{\ast}\|_{L^{p}(0,|\Omega|)}=\|f^{\#}\|_{L^{p}(\Omega^{\#})}\,;
\]
furthermore, the classical Hardy-Littlewood inequality holds true
\begin{equation}
\int_{\Omega}\vert f(x)\,  g(x)  \vert dx\leq\int_{0}^{\left\vert \Omega\right\vert}f^{\ast}(\sigma)\,  g^{\ast}(\sigma)  d\sigma=\int_{\Omega^{\#}}f^{\#}(x)\,g^{\#}(x)\,dx\,,
\label{HardyLit}%
\end{equation}
where $f,g$ are measurable functions on $\Omega$.

Here we recall an important tool in the proof of our main result, namely the following generalization of the \emph{Riesz rearrangement inequality} (see \cite[Theorem 2.2]{ALIEB} and \cite[Theorem 1]{Hichem} for a possible generalization).
\begin{theorem}
Let $F:\R^{+}\times\R^{+}\rightarrow\R^{+}$ be a continuous function such that $F(0,0)=0$ and
\begin{equation}
F(u_{2},v_{2})+F(u_{1},v_{1})\geq F(u_{2},v_{1})+F(u_{1},v_{2})\label{F}
\end{equation}
whenever $u_{2}\geq u_{1}>0$ and $v_{2}\geq v_{1}>0$.
Assume that $f, g$ are nonnegative measurable functions on $\R^{N}$ which satisfy \eqref{distribution}, then we have the inequalities
\begin{equation}
\int_{\R^{N}}\int_{\R^{N}}F(f(x),g(y))W(ax+by)\,dx\,dy\leq \int_{\R^{N}}\int_{\R^{N}}F(f^{\#}(x),g^{\#}(y))W(ax+by)\,dx\,dy\label{mainRieszineq}
\end{equation}
and
\[
\int_{\R^{N}}F(f(x),g(x))\,dx\leq \int_{\R^{N}}F(f^{\#}(x),g^{\#}(x))\,dx,
\]
for any nonnegative function $W\in L^{1}(\R^{N})$ and any choice of nonzero numbers $a$ and $b$.
\end{theorem}


\subsection{Mass concentration}
 The following definition based on mass concentration comparison will be widely used throughout the text. We refer the reader to \cite{ChRice}, \cite{ALTa}, \cite{Vsym82} for further details and related properties.
\begin{definition}
Let $f,g\in L^{1}_{loc}(\R^{N})$ be two radially symmetric functions on $\R^{N}$.
We say that $f$ is less concentrated than $g$, and we write
$f\prec
g$ if for
all $r>0$ we get
\[
\int_{B_{r}(0)}f(x)\,dx\leq \int_{B_{r}(0)}g(x)\,dx.
\]
\end{definition}
The partial order relationship $\prec$ is called \emph{comparison of mass concentrations}.
Of course, this definition can be suitably adapted if $f,g$ are defined in a ball $B_{R}(0)$ (considering the extension to zero outside $B_{R}(0)$). Moreover, if $f$ and $g$ are defined on two open sets of the same measure $\kappa$ we have that
$f^{\#}\prec g^{\#}$ if and only if
\[
\int_{0}^{\sigma}f^{\ast}(\tau)\,d\tau\leq \int_{0}^{\sigma}g^{\ast}(\tau)\,d\tau,
\]
for all $\sigma\in[0,\kappa]$.

The comparison of mass concentrations enjoys some nice equivalent formulations
(for the proof we refer to
\cite{Chong}, \cite{ALTa}, \cite{VANS05}).

\begin{lemma}\label{lemma1}
Let $f,g\in L_+^{1}(\Omega)$ two rearranged functions on a ball $\Omega=B_{R}(0).$ Then the following are equivalent:

\vskip7pt
\noindent(i) $f\prec g$;
\vskip7pt

\noindent(ii) for all $\phi\in L^\infty_+(\Omega)$,
$$\int_{\Omega}f(x)\phi(x)\,dx\leq \int_{\Omega^\#}f^\#(x)\phi^\#(x)\,dx.
$$
\vskip7pt

\noindent(iii) for all convex, nonnegative functions $\Phi:[0,\infty)\rightarrow [0,\infty)$ with $\Phi(0) = 0$ it holds
$$\int_{\Omega}\Phi(f(x))\,dx\leq \int_{\Omega}\Phi(g(x))\,dx.
$$
%
\end{lemma}
We explicitly observe that, if $f, g\in L^p(\Omega)$ $(1 < p \le\infty)$, then we may take $\phi \in L^{p'}(\Omega)$ in the point \textit{(ii)} above.

From this Lemma it easily follows that if $f$ and $g$ are $L^{1}$ functions on $\Omega$ such that $f^{\#}\prec g^{\#}$,
then
\begin{equation}
\|f\|_{L^{p}(\Omega)}\leq \|g\|_{L^{p}(\Omega)}\quad \forall p\in[1,\infty].
\end{equation}

\subsection{Functional spaces}\label{Functionspac}
Now we introduce the functional space where problem \eqref{mainpp} will be settled, namely a natural domain for the fractional $p$-Laplacian operator $(-\Delta _{p})^{s}$. Recall that $s\in (0,1)$ in our setting.  For any open set $\Omega$, we introduce the fractional Gagliardo seminorm
Gagliardo seminorm
\[
[u]_{W^{s,p}(\Omega)}=\left(\int_{\Omega}\int_{\Omega}\frac{|u(x)-u(y)|^{p}}{|x-y|^{N+sp}}dx\,dy\right)^{1/p},
\]
for a measurable function $u$ on $\Omega$. Then we define the fractional Sobolev space $W^{s.p}(\Omega)$ as the space
\[
W^{s,p}(\Omega)=\left\{u\in L^{p}(\Omega):\,[u]_{W^{s,p}(\Omega)}<\infty\right\},
\]
endowed with the norm
\[
\|u\|_{W^{s,p}(\Omega)}=\| u \|_{L^{p}(\Omega)}+[u]_{W^{s,p}(\Omega)}.
\]
We denote by $W_{0}^{s,p}(\Omega)$ the closure of $C_{c}^{\infty}(\Omega)$ in the $W^{s,p}(\Omega)$ topology.\\
The natural space when considering the operator $(-\Delta _{p})^{s}$ with homogeneous external Dirichlet condition will be denoted by $\widetilde{W}_{0}^{s,p}(\Omega)$, which is defined as
\[
\widetilde{W}_{0}^{s,p}(\Omega)=\left\{u:\R^{N}\rightarrow\R:\,[u]_{W^{s,p}(\R^{N})}<+\infty\text{ and }u=0 \text{ in }\R^{N}\setminus\Omega\right\}.
\]
When $p>1$ and $\Omega$ is an open bounded set with Lipschitz boundary, it can be proved that (see \cite[Proposition B.1]{BrasParSquas}) $\widetilde{W}_{0}^{s,p}(\Omega)$ coincides with the completion of $C_{c}^{\infty}(\Omega)$ with respect to the seminorm $ [\cdot]_{W^{s,p}(\R^{N})}$.  Moreover, when $sp\neq1$, it can also proved that $\widetilde{W}_{0}^{s,p}(\Omega)$ \emph{coincides} with $W_{0}^{s,p}(\Omega)$ (see \cite[Proposition B.1]{Brasco}), while in general for $sp=1$, we have a \emph{strict} inclusion $$\widetilde{W}_{0}^{s,p}(\Omega)\subset W_{0}^{s,p}(\Omega)$$ (see \cite[Remark 2.1]{secondeigenbrasco}).\\
A consequence of fractional Poincar\'e inequalities (see \cite[Lemma 2.4]{Brasco}) is that we can equip the space $\widetilde{W}_{0}^{s,p}(\Omega)$ with the Gagliardo seminorm
\[
\|u\|_{\widetilde{W}_{0}^{s,p}(\Omega)}=[u]_{W^{s,p}(\R^{N})}=\left(\int_{\R^{N}}\int_{\R^{N}}\frac{|u(x)-u(y)|^{p}}{|x-y|^{N+sp}}dx\,dy\right)^{1/p}.
\]

We finally recall the definition of Lorentz space (see, {\sl e.g.}, \cite{hunt}, \cite{oneil}). We introduce the \emph{maximal function} of Hardy and Littlewood associated to $v^{\ast}$, namely
\[
\bar{v}=\frac{1}{s}\int_{0}^{s}v^{\ast}(\sigma)d\sigma,
\]
and define for $1<p\leq\infty$, $0<q\le\infty$,
$$\|v\|_{p,q}=\left\{\begin{array}{ll}
\displaystyle\left(\int_0^{|\Omega|} \left(\bar{v}(s)\,s^{1 \over
p}\right)^q\,{ds\over
s}\right)^{1\over q},&
\qquad \text{if  } 0< q<\infty,\\
&\\
\displaystyle\sup_{s>0}\bar{v}(s)\,s^{1\over p},&
\qquad \text{if } q=\infty.
\end{array}
\right.
$$
The Lorentz space $L^{p,q}(\Omega)$ is defined as the set of the measurable
functions
$v$ such that
$\|v\|_{L^{p,q}(\Omega)}<+\infty.$ We will use the Lorentz spaces in the proof of Corollary \eqref{regularity}.

It is easy to verify that Lorentz spaces coincide with Lebesgue
spaces $L^{p}(\Omega)$ when $p=q$ and with Marcinkiewicz spaces $M^{p}(\Omega)$
when $q=\infty$. As regards the other values of the second index $q$, Lorentz spaces are
intermediate spaces between Lebesgue spaces
in the sense that, since $\Omega$ is bounded, the following inclusions hold true
\begin{align}\notag
&L^{p_2,r}(\Omega)\subset L^{p,q_1}(\Omega)\subset
L^{p,p}(\Omega)=L^{p}(\Omega) \subset\\
\notag
&\qquad \qquad \subset L^{p,q_2}(\Omega)\subset L^{p,\infty}(\Omega)=
M^{p}(\Omega)\subset L^{p_1,r}(\Omega),
\end{align}
when $1<p_1<p<p_2<\infty$,
$1\le q_1<p<q_2\le\infty$ and $1 \leq r\leq \infty$.

\subsection{Hypergeometric functions}\label{Hyperg}

We now recall the definition of the hypergeometric function $_{2}F_1(a,b;c;x)$ (see, for example, \cite[Ch. II]{magnus} fo further details). The definition of
$_{2}F_1(a,b;c;x)$ is given by
\begin{equation}\label{verydefHyp}
_{2}F_1(a,b;c;x)=\frac{\Gamma(c)}{\Gamma(b)\Gamma(a)}\sum \frac{\Gamma(a+n)\,\Gamma(b+n)}{\Gamma(c+n)}\frac{x^{n}}{n!}
\end{equation}
where the series converges in the unit interval $|x|<1$.
For $c>b>0$ and $0<\tau<1$, we have the following representation
\begin{equation}\label{represent}
_{2}F_1(a,b;c;x)=\frac{\Gamma(c)}{\Gamma(b)\Gamma(c-b)}\int_0^1\tau^{b-1}(1-\tau)^{c-b-1}(1-x\tau)^{-a}d\tau
\end{equation}
Some classical results about the derivatives of $_{2}F_1(a,b;c;x)$ read as
\begin{align*}\allowdisplaybreaks
_{2}F_1'(a,b;c;x)&=\frac {ab}c {\ }_{2}F_1(a+1,b+1;c+1;x),\displaybreak[1]\\ \\
_{2}F_1(a+1,b;c+1;x)&=\frac c{c-b} {\ }_{2}F_1(a,b;c;x)-\frac c{c-b}\frac{1-x}a{\ }_{2}F_1'(a,b;c;x),\displaybreak[1]\\ \\
_{2}F_1(a-1,b;c-1;x)&=\frac {c-1-bx}{c-1}{\ }_{2}F_1(a,b;c;x)+\frac{x(1-x)}{c-1}{\ }_{2}F_1'(a,b;c;x).
\end{align*}
A direct consequence of the above equalities is the following:
\begin{equation}\label{formulona}
_{2}F_1'(a,b;c;x)=\frac {ab}{c}{\ }_{2}F_1(a+1,b;c+1;x)+\frac {ax}{c}{\ }_{2}F'_1(a+1,b;c+1;x).
\end{equation}
and
\begin{equation}
\label{formulonaa}
x{\ }_{2}F_1'(a,b;c;x)+a{\ }_{2}F_1(a,b;c;x)=a{\ }_{2}F_1(a+1,b;c+1;x)+\frac {ax}{c}{\ }_{2}F'_1(a+1,b;c+1;x)
\end{equation}
In the present paper we have to use hypergeometric function in order to represent an integral which comes into play when one wants to calculate the fractional $p$-laplacian of radial function. Indeed,
 the following equality  holds true for $b>0$ and $|x|<1$ (see \cite[Ch. II, sub. 2.5.1]{magnus})
\begin{equation}\label{hyper}
\int_{0}^{\pi}\frac{\sin^{2b-1}\theta}{(1-2x\cos\theta+x^{2})^{a}}d\theta=\frac{\sqrt\pi\;\Gamma(b)}{\Gamma(b+\frac12)}{\ }_{2}F_1(a,a-b+\tfrac12;b+\tfrac12;x^{2}).
\end{equation}

Finally we recall that a direct computation in \eqref{verydefHyp} gives
\begin{equation}\label{zero}
{}_{2}F_{1}(a,b;c;0)=1,
\end{equation}
and, when $c>a+b$, for positive $a, b$, the following formula holds
 (see \cite[Ch. II, pag.40]{magnus})
\begin{equation}\label{gauss}
{}_{2}F_{1}(a,b;c;1)={\frac {\Gamma (c)\Gamma (c-a-b)}{\Gamma (c-a)\Gamma (c-b)}}.
\end{equation}
Using equality (see \cite[Ch. II, subs. 2.4.1]{magnus}), the last information allows to find
\begin{equation}\label{linear}
{}_{2}F_{1}(a,b;c;x)=(1-x)^{c-a-b}
{}_{2}F_{1}(c-a,c-b;c;x),
\end{equation}
which will be used in order to establish some asymptotic behaviours in the proof of Theorem \ref{maint}.

\section{Main results}\label{Main}
Assume that $s\in (0,1)$, $p\geq2$ and let $\Omega$ be a bounded open set with Lipschitz boundary.  As mentioned in the introduction, we will focus on the nonlinear Dirichlet problem
\begin{equation} \label{eq.0}
\left\{
\begin{array}[c]{lll}%
\left( -\Delta_{p}\right)^{s}u=f & & \text{in }%
\Omega,\\
\\
u=0 & & \text{on }\R^{N}\setminus\Omega.
\end{array}\right.
\end{equation}
A \emph{weak} solution to problem \eqref{eq.0} is a function $u\in \widetilde{W}_{0}^{s,p}(\Omega)$ such that
\begin{equation}\label{weak}
\frac{\gamma(N,s,p)}{2}\iint_{\R^{2N}}\frac{|u(x)-u(y)|^{p-2}(u(x)-u(y))(\varphi(x)-\varphi(y))}{|x-y|^{N+sp}}dx\,dy=\int_{\Omega}f(x)\varphi(x)dx,
\end{equation}
for all the test functions $\varphi\in \widetilde{W}_{0}^{s,p}(\Omega)$.
We will assume that the source term $f$ will possess enough summability in order to ensure that problem \eqref{eq.0} has a unique weak solution $u\in\widetilde{W}_{0}^{s,p}(\Omega)$. To this aim, we will assume that $f\in L^{m}(\Omega)$, where $m$ is such that
\begin{equation}\label{assumpt}
m\geq \frac{pN}{(p-1)N+sp} \text{ if }sp<N,\quad m>1 \text{ if }sp=N,\quad m\geq1 \text{ if }sp>N.
\end{equation}

Under the assumptions \eqref{assumpt} it is not difficult to show that the strictly convex functional
\[
\mathcal{J}(u)=\frac{\gamma(N,s,p)}{2p}\int_{\R^{N}}\int_{\R^{N}}\frac{|u(x)-u(y)|^{p}}{|x-y|^{N+sp}}dx\,dy-\int_{\Omega}fu\,dx.
\]
admits a minimizer in $\widetilde{W}_{0}^{s,p}(\Omega)$, which is in turn the unique weak solution to \eqref {eq.0}. Indeed, if $sp<N$, observe first that by the fractional Sobolev embedding (see \emph{e.g.} \cite[Theorem 6.5]{hitch} and \cite{Brasco})
\[
\widetilde{W}_{0}^{s,p}(\Omega)\hookrightarrow L^{p^{\ast}_{s}}(\Omega)
\]
where
\[
p^{\ast}_{s}=\frac{Np}{N-sp}
\]
and by the Young inequality,  one has for any $\varepsilon>0$ and for some positive constants $C$, $C^{\prime}$
\[
\int_{\Omega}f u\, dx\leq C\| f\|_{L^{(p^{\ast})^{\prime}}(\Omega)}^{p^{\prime}}+C^{\prime}\varepsilon[u]_{W^{s,p}(\R^{N})}^{p},
\]
hence for $\varepsilon$ small enough we have that $\mathcal{J}$ is bounded from below and that the following inequality holds for some positive constants $C$ and $C_{1}$:
\[
\mathcal{J}(u)\geq C_{1}\int_{\R^{N}}\int_{\R^{N}}\frac{|u(x)-u(y)|^{p}}{|x-y|^{N+sp}}dx\,dy-C\| f\|_{L^{(p^{\ast})^{\prime}}(\Omega)}^{p^{\prime}}.
\]
This last inequality is sufficient to show that from any minimizing sequence of $\mathcal{J}$ it is possible to extract a subsequence which is weakly converging to some $u$ in $\widetilde{W}_{0}^{s,p}(\Omega)$, that turns to be a minimizer of $\mathcal{J}$ by the lower semicontinuity of the norm in $\widetilde{W}_{0}^{s,p}(\Omega)$. A similar argument can be reproduced in the case $sp>N$, when we have the continuous embedding
\[
\widetilde{W}_{0}^{s,p}(\Omega)\hookrightarrow L^{\infty}(\Omega)\cap C^{0,\alpha}(\overline{\Omega})
\]
where $\alpha=1-N/sp$, see for instance \cite[Proposition 2.9]{Brasco}. This remark also applies to the case $sp=N$, where we have the continuous embedding
\[
\widetilde{W}_{0}^{s,p}(\Omega)\hookrightarrow L^{q}(\Omega)
\]
for all $q\in[\frac{N}{s},+\infty)$, see \cite[Theorem 6.9]{hitch} and \cite{Brasco}.
\\[0.6pt]
Our aim is to compare $u$ with the solution $v$ to a radial problem in the ball $\Omega^\#$. More precisely in the proof of our result an essential role will be played by the integral mean function of $u$ on balls, \emph{i.e},
\begin{equation}\label{sphmean}
U(x)=U(|x|)=\frac1{|x|^N}\int_0^{|x|}\mathfrc{u}(\rho)\rho^{N-1}d\rho,
\end{equation}
where $\mathfrc{u}(\rho)$ denotes the radial profile of $u^\#(x)$, that is, $\mathfrc{u}(|x|)=u^\#(x)$.

\begin{theorem}\label{maint}
Assume that $p\ge2$, $s\in(0,1)$ and $f\in L^{m}(\Omega)$ where $m$ satisfies one of the assumptions \eqref{assumpt}. Let $u\in  \widetilde{W}_{0}^{s,p}(\Omega)$ the weak solution to the Dirichlet problem \eqref{eq.0}.
Let $v\in  \widetilde{W}_{0}^{s,2}(\Omega^{\#})$ be the weak solution to the radial problem
\begin{equation} \label{eq.symmetric}
\left\{
\begin{array}[c]{lll}%
\left( -\Delta\right)^{s}v=g & & \text{in }%
\Omega^{\#},\\
\\
v=0 & & \text{on }\R^{N}\setminus\Omega^{\#},
\end{array}\right.
\end{equation}
where the datum $g=g(|x|)$ is the radial function defined by (we set as usual $r=|x|$)
\begin{equation}\label{g}
\begin{split}
g(r)=\mathsf{H}(N,s,p)\,r^{\frac{(N-s)(p-2)}{p-1}}&\left[\frac{(N-s)(p-2)}{p-1}\,\frac{1}{r^{N}}\left(\int_{B_r}f^{\#}dx\right)^{\frac{1}{p-1}}\right.\\
&\left.+\frac{N\omega_{N}}{p-1}\left(\int_{B_r}f^{\#}dx\right)^{\frac{2-p}{p-1}}f^{\#}(x)
\right],
\end{split}
\end{equation}
with
\[
\mathsf{H}(N,s,p)=\frac{\gamma(N,s,2)}{N\omega_{N}}\,\frac{\left(\mathcal{P}_{s}(B_{1})\right)^{\frac{p-2}{p-1}}}{\gamma(N,s,p)^{\frac{1}{p-1}}},
\]
being
\[
\mathcal{P}_{s}(B_{1})=\int_{B_r}\int_{B_r^c}\frac{1}{|x-y|^{N+s}}\,dx\,dy
\]
the fractional $s$-perimeter of the unit ball, see for instance \cite{CafVal}.
Then we have
\begin{equation}\label{massconc}
u^{\#}\prec v.
\end{equation}
%
%
%
\end{theorem}
\begin{remark}
Notice that when $p=2$ we have $g(x)=f^{\#}(x)$, so the symmetrized problem coincides with the one appearing in the \emph{linear} case $p=2$, hence Theorem \ref{maint} reduces to \cite[Theorem 3.1]{FeroneVolzone}. If $p>2$, it is not difficult to show that problem \eqref{eq.symmetric} admits a unique weak solution.
Indeed, using the simple estimate
$$f^\#(r)\le \frac 1{\omega_Nr^N}\int_{B_r}f^\#(x)\,dx,$$
we observe that
\begin{equation}\label{intermg}
g(r)\le  Cr^{-\frac{s(p-2)}{p-1}}\left(\frac 1{r^N}\int_{B_r}f^\#(x)\,dx\right)^{\frac1{p-1}}\end{equation}
where $C$ denotes a constant which can change from line to line. Assume first that $sp<N$.
An easy use of the H\"older inequality gives
\begin{equation}
g(r)\leq Cr^{-\frac{N+ms(p-2)}{m(p-1)}}\|f\|_{L^{m}(\Omega)}^{1/(p-1)}.
\end{equation}
 Thus
\[
\left(\int_{\Omega^{\#}}|g(r)|^{\frac{2N}{N+2s}}dx\right)^{\frac{N+2s}{N}}\leq C \|f\|_{L^{m}(\Omega)}^{1/(p-1)}\left(\int_{\Omega^{\#}}\frac{1}{r^{\ell}}dx\right)^{\frac{N+2s}{N}}
\]
where
\[
\ell=\frac{N+ms(p-2)}{m(p-1)}\frac{2N}{N+2s}
\]
and it is easy to show that $\ell<N$ if and only if
\begin{equation}\label{newest}
m>\frac{2N}{N(p-1)+2s}.
\end{equation}
Observe that
\[
\frac{Np}{(p-1)N+sp}>\frac{2N}{N(p-1)+2s}
\]
when $p^{2}-3p+2>0$, a condition which is assured by the restriction $p>2$. Therefore \eqref{assumpt} yields \eqref{newest}.\\
As regards to the range $sp\geq N$, with $N\geq2$, or $N=1$ and $s\in(0,1/2)$, by \eqref{intermg} we get the inequality
\begin{equation}\label{boundg}
g(r)\leq Cr^{-\frac{N+s(p-2)}{(p-1)}}\|f\|_{L^{1}(\Omega)}^{1/(p-1)}
\end{equation}
which implies
\[
\left(\int_{\Omega^{\#}}|g(r)|^{\frac{2N}{N+2s}}dx\right)^{\frac{N+2s}{N}}\leq C \|f\|_{L^{1}(\Omega)}^{1/(p-1)}\left(\int_{\Omega^{\#}}\frac{1}{r^{\ell}}dx\right)^{\frac{N+2s}{N}}
\]
where
\[
\ell=\frac{N+s(p-2)}{p-1}\frac{2N}{N+2s}
\]
hence $\ell<N$ if and only if
\begin{equation}
\label{boundpN}
p>\frac{3N-2s}{N}
\end{equation}
is satisfied. But an easy computation shows that under our assumption $N>2s$ the following bound holds
\[
\frac{N}{s}>\frac{3N-2s}{N}
\]
thus \eqref{boundpN} follows. In all these cases we showed that $g\in L^{2N/(N+2s)}(\Omega^{\#})$, thus  the linear problem \eqref{eq.symmetric} has a unique solution. In the case $N=1$, $s\in [1/2,1)$ it will be sufficient to show that $g\in L^{q}(\Omega^{\#})$ for some $q>1$: this can be done by choosing
\[
q<\frac{p-1}{1+s(p-2)}
\]
in order to have
\[
\left(\int_{\Omega^{\#}}|g(r)|^{q}dx\right)^{\frac{1}{q}}\leq C \|f\|_{L^{1}(\Omega)}^{1/(p-1)}\left(\int_{\Omega^{\#}}\frac{1}{r^{\ell}}dx\right)^{\frac{1}{q}}<\infty
\]
where
\[
\ell=\frac{1+s(p-2)}{p-1}q<1.
\]

\end{remark}

As a consequence of Theorem \ref{maint} we will prove the following regularity result for solutions $u$ to \eqref{eq.0}.
\begin{corollary}\label{regularity}
Let us choose $s\in(0,1)$, $p\geq2$ and assume that $sp<N$. Let $u\in \widetilde{W}_{0}^{s,p}(\Omega)$ the weak solution to the Dirichlet problem \eqref{eq.0}.
We have:

\noindent 1. if $f\in L^{{m,\frac{Nm}{N+sm(p-2)}}}(\Omega)$, with $\frac{pN}{((p-1)N+sp)}\le m<N/(sp)$, then $u\in L^{q}(\Omega)$, with
\begin{equation}
q=\frac{Nm(p-1)}{N-smp}\label{q}
\end{equation}
and there exists a constant $C$ such that:
\[
\|u\|_{L^{q}(\Omega)}\leq C\|f\|_{L^{{m,\frac{Nm}{N+sm(p-2)}}}(\Omega)}^{\frac1{p-1}};
\]

\noindent 2.  if $f\in L^{{m}}(\Omega)$, with $m>N/(sp)$, then $u\in L^{\infty}(\Omega)$ and there exists a constant $C$ such that:
\[
\|u\|_{L^{\infty}(\Omega)}\leq C\|f\|_{L^{m}(\Omega)}^{\frac1{p-1}}.
\]
\end{corollary}
\begin{remark}
The case $m>N/(sp)$ in Corollary \ref{regularity} is obtained in \cite[Theorem 3.1]{secondeigenbrasco} and in \cite[Theorem 3.1]{BarPer} by Moser's iteration techniques . The case $\frac{pN}{((p-1)N+sp)}\le m<N/(sp)$ and $f\in L^{m}(\Omega)$ is studied in  \cite[Theorem 3.4]{BarPer}. Notice that the Lorentz space in case 1. of Corollary \ref{regularity} is smaller that $L^{m}(\Omega)$, thus our regularity result is not sharp: this discrepancy is due to the choice of the symmetrized problem \eqref{eq.symmetric}, which is linear.
\end{remark}

\section{Proofs}\label{Proofs}

First of all we prove a preliminary result which will be used in order to apply Riesz rearrangement inequality \eqref{mainRieszineq}. Let $\mathcal{G}_{t,h}$ , $t,h>0$, be the classical truncation function
\begin{equation}
\mathcal{G}_{t,h}(\theta)  =\left\{
\begin{array}
[c]{lll}%
h &  & \text{if }\theta > t+h\\
&  & \\
\theta-t\, &  & \text{if }t< \theta \le t+h\\
&  & \\
0 &  & \text{if }\theta \leq t.\text{ }%
\end{array}
\right.\label{truncation}
\end{equation}
The following Lemma establishes a property of a nonlinear function $F(u,v)$ involving the truncations $\mathcal{G}_{t,h}$, that will be employed as an application of inequality \eqref{mainRieszineq}, a basic ingredient in the proof of Theorem \ref{maint}.
\begin{lemma}\label{effe}
For $p\ge2$, let $F:\R^{+}\times\R^{+}\rightarrow\R$ be defined as
\begin{equation}\label{effedef}
F(u,v)=u^{p}+v^{p}-|u-v|^{p-2}(u-v)(\mathcal{G}_{t,h}(u)-\mathcal{G}_{t,h}(v)).
\end{equation}
Then $F(u,v)$ is a continuous function, with $F(0,0)=0$, such that
\begin{equation}\label{effepos}
F(u,v)\geq0
\end{equation}
and
\begin{equation}\label{effeprop}
F(u_{2},v_{2})+F(u_{1},v_{1})\geq F(u_{2},v_{1})+F(u_{1},v_{2})
\end{equation}
\end{lemma}

\noindent{\sc Proof.}
As regards \eqref{effepos} we observe that
\[
(u-v)(\mathcal{G}_{t,h}(u)-\mathcal{G}_{t,h}(v))\le (u-v)^2
\]
Indeed

\noindent-- if $0\le u\le t$ and $0\le v\le t$ then
$$(u-v)(\mathcal{G}_{t,h}(u)-\mathcal{G}_{t,h}(v))=(u-v)\cdot0\le (u-v)^2$$

\noindent-- if $0\le u\le t$ and $t<v<t+h$ then
$$(u-v)(\mathcal{G}_{t,h}(u)-\mathcal{G}_{t,h}(v))=(u-v)(t-v)\le (u-v)^2$$

\noindent-- if $0\le u\le t$ and $t+h\le v$ then
$$(u-v)(\mathcal{G}_{t,h}(u)-\mathcal{G}_{t,h}(v))=(u-v)(-h)\le (u-v)^2$$

\noindent-- if $t<u<t+h$ and $0\le v\le t$ then
$$(u-v)(\mathcal{G}_{t,h}(u)-\mathcal{G}_{t,h}(v))=(u-v)(u-t)\le (u-v)^2$$

\noindent-- if $t<u<t+h$ and $t<v<t+h$ then
$$(u-v)(\mathcal{G}_{t,h}(u)-\mathcal{G}_{t,h}(v))=(u-v)^2$$

\noindent-- if $t<u<t+h$ and $t+h\le v$ then
$$(u-v)(\mathcal{G}_{t,h}(u)-\mathcal{G}_{t,h}(v))=(u-v)(u-t-h)\le (u-v)^2$$

\noindent-- if $t+h\le u$ and $0\le v\le t$ then
$$(u-v)(\mathcal{G}_{t,h}(u)-\mathcal{G}_{t,h}(v))=(u-v)h\le (u-v)^2$$

\noindent-- if $t+h\le u$ and $t<v<t+h$ then
$$(u-v)(\mathcal{G}_{t,h}(u)-\mathcal{G}_{t,h}(v))=(u-v)(t+h-v)\le (u-v)^2$$

\noindent-- if $t+h\le u$ and $t+h\le v$ then
$$(u-v)(\mathcal{G}_{t,h}(u)-\mathcal{G}_{t,h}(v))=(u-v)\cdot0\le (u-v)^2$$

\noindent Then, by monotonicity,
\[
F(u,v)\ge u^{p}+v^{p}-|u-v|^{p}\geq0.
\]

In order to prove \eqref{effeprop} we observe that, for every fixed $u\ge0$ the function
\begin{equation*}
\Phi(v)=|u-v|^{p-2}\big(\mathcal{G}_{t,h}(u)-\mathcal{G}_{t,h}(v)\big),\quad v\ge0,
\end{equation*}
is decreasing with respect to $v$. Indeed, we can compute the derivative for a.e. $v$ to get
\begin{equation*}
\Phi'(v)=-|u-v|^{p-4}\big((p-2)(u-v)\big(\mathcal{G}_{t,h}(u)-\mathcal{G}_{t,h}(v)\big)+|u-v|^2\mathcal{G}'_{t,h}(v)\big)\le0.
\end{equation*}
This means that, for every $0\le v_1\le v_2$, it holds
\begin{equation*}
\Phi(v_1)\ge\Phi(v_2),
\end{equation*}
that is, for every $u\ge0$ and $0\le v_1\le v_2$
\begin{equation}\label{phi}
|u-v_1|^{p-2}\big(\mathcal{G}_{t,h}(u)-\mathcal{G}_{t,h}(v_1)\big)-
|u-v_2|^{p-2}\big(\mathcal{G}_{t,h}(u)-\mathcal{G}_{t,h}(v_2)\big)\ge0.
\end{equation}
Now, for every fixed $0\le v_1\le v_2$, we consider the function
\begin{equation*}
\Psi(u)=F(u,v_2)-F(u,v_1),\quad u\ge0.
\end{equation*}
Using the definition \eqref{effedef} of $F(u,v)$, we can compute the derivative of $\Psi(u)$ for a.e. $u$ to get
\begin{align*}
\Psi'(v)=&(p-1)\big(|u-v_1|^{p-2}\big(\mathcal{G}_{t,h}(u)-\mathcal{G}_{t,h}(v_1)\big)-
|u-v_2|^{p-2}\big(\mathcal{G}_{t,h}(u)-\mathcal{G}_{t,h}(v_2)\big)\\
\\
&\>+\mathcal{G}'_{t,h}(u)\big(|u-v_1|^{p-2}(u-v_1)-|u-v_2|^{p-2}(u-v_2)\big)\ge0,
\end{align*}
where we have used inequality \eqref{phi} and the monotonicity of the function $\phi(t)=|t|^{p-2}t$, $t\in\R$.

Using the monotonicity of $\Psi(u)$, we have, for $0\le u_1\le u_2$ and $0\le v_1\le v_2$,
\begin{equation*}
F(u_2,v_2)-F(u_2,v_1)\ge F(u_1,v_2)-F(u_1,v_1),
\end{equation*}
that is, \eqref{effeprop}.
\hfill$\square$
\vskip.3cm
\noindent{\sc Proof of Theorem \ref{maint}.}
We first consider the case $f\ge0$ and we assume $f\in C_{0}^{\infty}(\Omega)$. By \cite[Theorem 1.4]{brascohol} or \cite[Theorem 1.1]{IanMoscSquas} we have that the solution $u$ to \eqref{mainpp} is locally H\"older continuous in $\Omega$ (see also \cite[Corollary 1.1]{KuusiMing} where continuity is derived for more general nonlocal operators). If $\mathcal{G}_{t,h}$ is the truncation function in \eqref{truncation}, we use the test function
\[
\varphi(x)=\mathcal{G}_{t,h}(u(x)),
\]
in the weak formulation of \eqref{eq.0}, obtaining
\begin{align}
\displaystyle\frac{\gamma(N,s,p)}{2}&\int_{\R^{N}}\int_{\R^{N}}\frac{|u(x)-u(y)|^{p-2}\left(u(x)-u(y)\right)\left(\mathcal{G}_{t,h}(u(x))-\mathcal{G}_{t,h}(u(y))\right)}{|x-y|^{N+sp}}dx\,dy\label{eqtest}
\\
&=\int_{\Omega}f(x)\,\mathcal{G}_{t,h}(u(x))\,dx.\notag
\end{align}
We will prove the inequality
\begin{align}\label{Polyatype}
&\int_{\R^{N}}\int_{\R^{N}}\frac{|u(x)-u(y)|^{p-2}\left(u(x)-u(y)\right)\left(\mathcal{G}_{t,h}(u(x))-\mathcal{G}_{t,h}(u(y))\right)}{|x-y|^{N+sp}}dxdy
\\ &\geq
 \int_{\R^{N}}\int_{\R^{N}}\frac{|u^{\#}(x)-u^{\#}(y)|^{p-2}\left(u^{\#}(x)-u^{\#}(y)\right)\left(\mathcal{G}_{t,h}(u^{\#}(x))-\mathcal{G}_{t,h}(u^{\#}(y))\right)}{|x-y|^{N+sp}}dxdy.\nonumber
\end{align}
Following \cite[Section 9]{ALIEB}, we write
\begin{align*}
&\int_{\R^{N}}\int_{\R^{N}}\frac{|u(x)-u(y)|^{p-2}\left(u(x)-u(y)\right)\left(\mathcal{G}_{t,h}(u(x))-\mathcal{G}_{t,h}(u(y))\right)}{|x-y|^{N+sp}}dxdy\\
&=\frac{1}{\Gamma(\frac{N+sp}{2})}\int_{0}^{\infty}I_{\alpha}[u,t,h]\,\alpha^{(N+sp)/2-1}d\alpha,
\end{align*}
where
\begin{equation}
I_{\alpha}[u,t,h]=\int_{\R^{N}}\int_{\R^{N}}|u(x)-u(y)|^{p-2}\bigl(u(x)-u(y)\bigr)\bigl(\mathcal{G}_{t,h}(u(x))-\mathcal{G}_{t,h}(u(y))\bigr)\exp[-|x-y|^{2}\alpha]dx\,dy.
\label{applRiesz}
\end{equation}
We want to prove\begin{equation}
I_{\alpha}[u,t,h]\geq I_{\alpha}[u^{\#},t,h]\label{IneqI},
\end{equation}
for all $\alpha>0$. To this aim, we define the function $F(u,v)$ according to \eqref{effedef}. Then Lemma \ref{effe} ensures that $F$ is eligible in Riesz rearrangement inequality \eqref{mainRieszineq}, with the choice $W_{\alpha}(x)=\exp[-|x|^{2}\alpha]$ and $a=1,\,b=-1$.
Plugging such function $F$ in \eqref{mainRieszineq} yields
\[
\int_{\R^{N}}\int_{\R^{N}}F(u(x),u(y))\,W_{\alpha}(x-y)dx\,dy
\leq \int_{\R^{N}}\int_{\R^{N}}F(u^{\#}(x),u^{\#}(y))\,W_{\alpha}(x-y)dx\,dy,
\]
then the usual equimeasurability property of rearrangements and the symmetry of the kernel $W_{\alpha}$ applied to \eqref{applRiesz} gives \eqref{IneqI}.\\
As a consequence we have
\begin{align}
\frac{\gamma(N,s,p)}{2}&\int_{\R^{N}}\int_{\R^{N}}\frac{|u^{\#}(x)-u^{\#}(y)|^{p-2}\bigl(u^{\#}(x)-u^{\#}(y)\bigr)\bigl(\mathcal{G}_{t,h}(u^{\#}(x))-\mathcal{G}_{t,h}(u^{\#}(y))\bigr)}{|x-y|^{N+sp}}\,dy\,dx\label{mainineq}\\
 &\leq \int_{\Omega}f(x)\,\mathcal{G}_{t,h}(u(x))\,dx.\notag
 \end{align}
Now we set
 \[
 \mathfrc{u}(x)=\mathfrc{u}(|x|):=u^{\#}(x),
 \]
hence $\mathfrc{u}$ is a nonincreasing continuous function defined on $\R^N$ vanishing for $|x|\ge R>0$, $R$ being the radius of $\Omega^{\#}$. For any $0\le t\le u_{\text{max}}$ there exists a unique $r(t)$ such that $|\{x:\mathfrc{u}(x)>t\}|=|B_{r(t)}(0)|$.
{
Then we check how to pass to the limit as $h\rightarrow 0$ in \eqref{mainineq}. Let us consider the following integral:
\begin{align}\allowdisplaybreaks
I_{t,h}= \frac1{N\omega_{N}}
\int_{\R^{N}}\int_{\R^{N}}\frac{|\mathfrc{u}(x)-\mathfrc{u}(y)|^{p-2}\left(\mathfrc{u}(x)-\mathfrc{u}(y)\right)\left(\mathcal{G}_{t,h}(\mathfrc{u}(x))-\mathcal{G}_{t,h}(\mathfrc{u}(y))\right)}{|x-y|^{N+sp}}dxdy,
\label{integralh}
\end{align}
thus \eqref{mainineq} can be rewritten as
\begin{equation}
 N\omega_{N}\frac{\gamma(N,s,p)}{2} I_{t,h}\leq \int_{\Omega}f(x)\,\mathcal{G}_{t,h}(u(x))\,dx.\label{maininequalith}
\end{equation}

Putting  $r=|x|$, we have:
\begin{equation*}
\begin{split}
I_{t,h}=\int_0^{+\infty}&\left(\int_0^{+\infty}|\mathfrc{u}(r)-\mathfrc{u}(\rho)|^{p-2}
\big(\mathfrc{u}(r)-\mathfrc{u}(\rho)\big)\right.\\
&\times\left.\big(\mathcal{G}_{t,h}(\mathfrc{u}(r))-\mathcal{G}_{t,h}(\mathfrc{u}(\rho))\big)\Theta_{N,s,p}(r,\rho)\rho^{N-1}d\rho\right)r^{N-1}dr,
\end{split}
\end{equation*}
where
\begin{equation}\label{phi}
\Theta_{N,s,p}(r,\rho)=\frac1{N\omega_{N}}\int_{|x'|=1}\left(\int_{|y'|=1}\frac1{|r\,x'-\rho\,y'|^{N+sp}}dH^{N-1}(y')\right)dH^{N-1}(x').
\end{equation}
Thus from \eqref{hyper} (see \cite{FeroneVolzone}) it follows
\begin{equation}\label{explicit}
\Theta_{N,s,p}(r,\rho)=
\left\{
\begin{array}{ll}
\dfrac{\alpha_{N}}{\rho^{N+sp}}\>{}_{2}F_1\left(\dfrac{N+sp}2,\dfrac{sp}{2}+1;\dfrac N2;\dfrac{r^{2}}{\rho^{2}}\right)  &\quad   \text{if }0\le r<\rho<+\infty \\
&\\
\dfrac{\alpha_{N}}{r^{N+sp}}\>{}_{2}F_1\left(\dfrac{N+sp}2,\dfrac{sp}{2}+1;\dfrac N2;\dfrac{\rho^{2}}{r^{2}}\right) &   \quad   \text{if }0\le \rho<r<+\infty,\\
\end{array}
\right.
\end{equation}
where
$$\alpha_{N}=\frac{2\pi^{\frac{N-1}2}}{\Gamma(\frac{N-1}2)}.
$$
Moreover by \eqref{zero} we have the following asymptotic behaviors
\begin{equation}\label{infty}
\left\{
\begin{array}{ll}
\Theta_{N,s,p}(r,\rho)\sim\dfrac{\alpha_{N}}{r^{N+sp}}&\qquad\text{ as }r\rightarrow+\infty\\
\\
\Theta_{N,s,p}(r,\rho)\sim\dfrac{\alpha_{N}}{\rho^{N+sp}}&\qquad\text{ as }\rho\rightarrow+\infty,
\end{array}
\right.
\end{equation}
and a combination of \eqref{linear} and \eqref{gauss} provides
\begin{equation}\label{uno}
\Theta_{N,s,p}(r,\rho)\sim\dfrac1{|r-\rho|^{1+sp}}\qquad\text{ as }|r-\rho|\rightarrow0.
\end{equation}
Then, since $\mathfrc{u}$ is radially decreasing, all the machinery contained in the previous paper \cite{FeroneVolzone} can be used to employ the Lebesgue monotone convergence in order to pass to the limit as $h\rightarrow0$ in \eqref{maininequalith} and find the inequality
\begin{equation}\label{inequality}
\gamma(N,s,p)
\int_0^{r}\left(\int_{r}^{+\infty}|\mathfrc{u}(\tau)-\mathfrc{u}(\rho)|^{p-1}\Theta_{N,s,p}(\tau,\rho)\rho^{N-1}d\rho\right)\tau^{N-1}d\tau
\le\int_0^{r}f^{*}(\omega_{N}\rho^{N})\rho^{N-1}d\rho.
\end{equation}
Observe that by \eqref{maininequalith} the quotient ratio $I_{t,h}/h$ remains bounded, therefore the integral at the left hand side of \eqref{inequality} is finite. Now, writing \eqref{inequality} in cartesian coordinates again,
\begin{equation}
\gamma(N,s,p)\int_{B_r}\int_{B_r^c}\frac{|u^{\#}(x)-u^{\#}(y)|^{p-1}}{|x-y|^{N+sp}}\,dx\,dy
 \leq \int_{B_r}f^\#(x)\,dx.\label{mainineqnon}
 \end{equation}
At this point we write
\[
N+2s=\frac{N+2s}{p-1}+(N+s)\frac{p-2}{p-1}
\]
and use H\"older inequality in \eqref{mainineqnon} with exponents $1/(p-1)$ and $(p-2)/(p-1)$ to obtain
\begin{align}\label{newnon}
&\int_{B_r}\int_{B_r^c}\frac{|u^{\#}(x)-u^{\#}(y)|}{|x-y|^{N+2s}}\,dx\,dy\nonumber \\ & \le \left(\int_{B_r}\int_{B_r^c}\frac{|u^{\#}(x)-u^{\#}(y)|^{p-1}}{|x-y|^{N+sp}}\,dx\,dy\right)^{\frac1{p-1}}
\left(\int_{B_r}\int_{B_r^c}\frac{1}{|x-y|^{N+s}}\,dx\,dy\right)^{\frac{p-2}{p-1}}\nonumber\\&
 =r^{\frac{(N-s)(p-2)}{p-1}}\mathcal{P}_{s}(B_{1})^{\frac{p-2}{p-1}}\,\left(\int_{B_r}\int_{B_r^c}\frac{|u^{\#}(x)-u^{\#}(y)|^{p-1}}{|x-y|^{N+sp}}\,dx\,dy\right)^{\frac1{p-1}},
 \end{align}
 being
$\mathcal{P}_{s}(B_{1})$ the fractional perimeter of the unit ball. Then by \eqref{mainineqnon} we have
\begin{equation}\label{intermineq}
\int_{B_r}\int_{B_r^c}\frac{|u^{\#}(x)-u^{\#}(y)|}{|x-y|^{N+2s}}\,dx\,dy\leq
r^{\frac{(N-s)(p-2)}{p-1}}\frac{\left(\mathcal{P}_{s}(B_{1})\right)^{\frac{p-2}{p-1}}}{\gamma(N,s,p)^{\frac{1}{p-1}}}
\left(
\int_{B_{r}}f^{\#}\,dx\right)^{\frac{1}{p-1}}.
\end{equation}
Now, arguing as in \cite{FeroneVolzone} we find that
\[
\int_{B_r}\int_{B_r^c}\frac{|u^{\#}(x)-u^{\#}(y)|}{|x-y|^{N+2s}}\,dx\,dy=\frac{N\omega_{N}}{\gamma(N,s,2)} r^{N} (-\Delta)_{\R^{N+2}}^{s}U(r),
\]
where $(-\Delta)_{\R^{N+2}}^{s}$ denotes the $s$-laplacian computed on a radial function in $\R^{N+2}$.
Thus \eqref{intermineq} yields
 \begin{align*}
(-\Delta)_{\R^{N+2}}^{s}U(r)&\leq\mathsf{H}(N,s,p)\, r^{\frac{(N-s)(p-2)}{p-1}}\frac{1}{r^N}\left(\int_{B_r}f^\#(x)\,dx\right)^{\frac1{p-1}} \nonumber\\& =\mathsf{H}(N,s,p)\,\frac{1}{r^{\frac N{p-1}+s\frac{p-2}{p-1}}}\left(\int_{B_r}f^\#(x)\,dx\right)^{\frac1{p-1}}. \label{comparison}
\end{align*}


Now we observe that a direct computation (see the proof of \cite[Theorem 3.1]{FeroneVolzone}) shows that the solution $v$ to problem \eqref{eq.symmetric} is such that the integral mean function of $v$
\begin{equation}\label{sphmeanv}
V(x)=V(|x|)=\frac1{|x|^N}\int_0^{|x|}v(\rho)\rho^{N-1}d\rho,
\end{equation}
satisfies
\begin{align*}
(-\Delta)_{\R^{N+2}}^{s}V(x)&=\frac1{|x|^N}\int_0^{|x|}g(\rho)\rho^{N-1}\,d\rho\\
&=\mathsf{H}(N,s,p)\frac1{|x|^N}\int_0^{|x|} \frac{d}{d\rho}\left[r^{\frac{(N-s)(p-2)}{p-1}}\left(\int_{B_{\rho}}f^{\#}(x)dx\right)^{\frac{1}{p-1}}\right]\,d\rho\nonumber\\& =\mathsf{H}(N,s,p)\,\frac{1}{|x|^{\frac N{p-1}+s\frac{p-2}{p-1}}}\left(\int_{B_r}f^\#(x)\,dx\right)^{\frac1{p-1}}\nonumber\\&
=\frac{1}{|x|^{s\frac{p-2}{p-1}}}\left(\frac{1}{|x|^{N}}\int_{B_r}f^\#(x)\,dx\right)^{\frac1{p-1}},
\label{comparison}
\end{align*}
which provides the radially decreasing monotonicity of $V$ in $\R^{N+2}$.
It follows that
\[
(-\Delta)_{\R^{N+2}}^{s}U(r)\leq (-\Delta)_{\R^{N+2}}^{s}V(r)
\]
and we can apply the comparison principle for the fractional Laplacian (see \cite[Theorem 3.1]{FeroneVolzone} again) , which gives
\[
U\leq V,
\]
namely \eqref{massconc}.

The result is then achieved when $f\ge0$ and $f$ is regular. It is possible to remove the regularity assumption by using a suitable sequence of data as made, for example in \cite[Section 5.2]{FeroneVolzone}.

As regards the sign assumption, one can observe that
the comparison principle (see, for example, \cite{BarPer}) states that $|u|\le\tilde u$, being $\tilde u$ the solution to problem \eqref{mainpp} having $|f|$ as source term. Thus, we have:
$$u^\#\prec \tilde u^\#\prec v$$
and the theorem is completely proved.
\hfill$\square$
\vskip.3cm

\begin{remark}
Now we make a closer inspection to the Holder inequality \eqref{newnon}, which we can write in the more representative form (we set $r=1$)
\begin{equation}\label{nolocH}
\begin{split}
&\int_{B_{1}}\left(P.V.\int_{\R^{N}}\frac{u(x)-u(y)}{|x-y|^{N+2s}}dy\right)\,dx\\ &\leq \left(P_{s}(B_{1})\right)^{\frac{p-2}{p-1}}\left(\int_{B_{1}}\left(P.V.\int_{\R^{N}}\frac{|u(x)-u(y)|^{p-2}(u(x)-u(y))}{|x-y|^{N+ps}}dy\right)\,dx\right)^{\frac{1}{p-1}}.
\end{split}
\end{equation}
This inequality, is in general strict for $s<1$. We want to use an heuristic argument to show that this equality will not degenerate to an \emph{equality} when $s\rightarrow1$.
\\
To this aim, we check the asymptotics as $s\rightarrow1$ of \eqref{nolocH}, taking into account that, for sufficiently smooth functions,
\[
\lim_{s\rightarrow1^{-}}(1-s)\,P.V.\int_{\R^{N}}\frac{u(x)-u(y)}{|x-y|^{N+2s}}dy=\frac{\pi^{N/2}}{4\Gamma\left(\frac{N+2}{2}\right)}(-\Delta)u(x),
\]
and ( see \cite[Lemma 5.1]{CDVAZ})
\[
\lim_{s\rightarrow1^{-}}(1-s)\,P.V.\int_{\R^{N}}\frac{|u(x)-u(y)|^{p-2}(u(x)-u(y))}{|x-y|^{N+ps}}dy=\frac{\pi^{\frac{N-1}{2}}\,\Gamma\left(\frac{p+1}{2}\right)}
{p\Gamma\left(\frac{N+p}{2}\right)}(-\Delta_{p})u(x),
\]
and finally by \cite[Theorem 4]{Ludwig}
\[
\lim_{s\rightarrow1^{-}}(1-s)\mathcal{P}_{s}(B_{1})=\omega_{N-1}\mathcal{P}(B_{1}).
\]
Passing to the limit as $s\rightarrow 1$ in \eqref{nolocH}  we find
\[
\frac{\pi^{N/2}}{4\Gamma\left(\frac{N+2}{2}\right)}\int_{B_{1}}(-\Delta)u(x)\,dx\leq
(\mathcal{P}(B_{1}))^{\frac{p-2}{p-1}}\,\frac{\pi^{\frac{N-1}{2}}}{\Gamma\left(\frac{N+1}{2}\right)^{\frac{p-2}{p-1}}}
\left(\frac{\Gamma\left(\frac{p+1}{2}\right)}{p\Gamma\left(\frac{N+p}{2}\right)}\right)^{\frac{1}{p-1}}\left(\int_{B_{1}}(-\Delta_{p})u\,dx\right)^{\frac{1}{p-1}}
\]
But now we observe that a formal use of the divergence theorem and the radiality of $u$ yields
\[
\int_{B_{1}}(-\Delta)u\,dx=(\mathcal{P}(B_{1}))^{\frac{p-2}{p-1}}\left(\int_{B_{1}}(-\Delta_{p})u\,dx\right)^{\frac{1}{p-1}},
\]
but the presence of the constants depending on $p$ in the previous inequality prevents \eqref{nolocH} to degenerate into an equality in the limit $s\rightarrow1$.\\
On the other hand, we observe that \eqref{mainineqnon} can be rewritten as
\[
\int_{B_{r}}(-\Delta_p)^{s}u\,dx\leq \int_{B_{r}}f^{\#}dx,
\]
then formally passing to the limit as $s\rightarrow1$ gives
\[
\int_{B_{r}}(-\Delta_{p})u\,dx\leq \int_{B_{r}}f^{\#}dx,
\]
from which
\[
-\frac{u^{\prime}(r)}{r}\leq \frac{1}{r^{\frac{N-1}{p-1}+1}}\left(\int_{B_{r}}f^{\#}\right)^{\frac{1}{p-1}}=-\frac{v^{\prime}(r)}{r}
\]
that is the classical Talenti's result in \cite{talNON}. The previous remark then shows that using H\"{o}lder does not allow to recover the same inequality in the asymptotic limit $s\rightarrow1$.
\end{remark}

\noindent{\sc Proof of Corollary \ref{regularity}.}
Recall that $g$ satisfies the bound \eqref{boundg}.

If $f\in L^{{m,\frac{Nm}{N+sm(p-2)}}}(\Omega)$, with $\frac{pN}{((p-1)N+sp)}\le m<N/(sp)$, and $t=\frac{Nm(p-1)}{N+sm(p-2)}$, we have:
\begin{equation}\label{norms}
\|g\|_{L^t(\Omega)}\le C\left(\int_0^{+\infty}\tau^{-\frac{st(p-2)}{N(p-1)}}\left(\frac 1{\tau}\int_0^\tau f^*(\sigma)\,d\sigma\right)^{\frac t{p-1}} d\tau\right)^{\frac1t}=C\|f\|_{L^{m,\frac{Nm}{N+sm(p-2)}}(\Omega)}^{\frac1{p-1}}.
\end{equation}
Observe that $t\geq2N/(N+2s)$ exactly when \eqref{newest} holds.} This guarantees that the linear problem \eqref{eq.symmetric} has a unique solution.
By \cite[Theorem 3.2]{FeroneVolzone} we know that
\begin{equation}\label{LpreguFerVol}
\|v\|_{L^q(\Omega)}\le C\|g\|_{L^t(\Omega)}, \quad\text{with }q=\frac{Nt}{N-2st}=\frac{Nm(p-1)}{N-smp}.
\end{equation}
Now by Hardy-Littlewood inequality \eqref{HardyLit} we have $v\prec v^{\#}$, therefore from Theorem \ref{maint} we easily infer
\[
u^{\#}\prec v^{\#}
\]
and then \eqref{LpreguFerVol} and Lemma \ref{lemma1} imply
$$\|u\|_{L^q(\Omega)}\le\|v\|_{L^q(\Omega)}\le \|f\|_{L^{{m,\frac{Nm}{N+sm(p-2)}}}(\Omega)}^{\frac1{p-1}}.
$$

If $f\in L^{{m}}(\Omega)$, with $m>N/(sp)$, we choose $m'$ with $N/(sp)<m'<m$ and we put $t=\frac{Nm'(p-1)}{N+sm'(p-2)}$. It follows that
$$t>\frac N{2s},$$
so, by \cite[Theorem 3.2]{FeroneVolzone} we have
$$\|u\|_{L^\infty(\Omega)}\le C\|g\|_t$$
and then by \eqref{norms}
$$\|u\|_{L^\infty(\Omega)}\le\|v\|_{L^\infty(\Omega)}\le C \|f\|_{L^{{m',\frac{Nm'}{N+sm'(p-2)}}}(\Omega)}^{\frac1{p-1}}\le C\|f\|_{L^m(\Omega)}^{\frac1{p-1}},
$$
where the last inequality comes from the well-known inclusions in Lorentz spaces when $m'<m$, see Section \ref{Functionspac}.

\hfill$\square$
\section{Comments and open problems}\label{Open problems}
\noindent {$\bullet$} As it was mentioned in the introduction, the degenerate condition $p>2$ seems to be indispensable for the validity of Lemma \eqref{effe}, where a basic importance is given to such convexity property. Therefore, it would be extremely interesting to have an extension of a Theorem in the form \eqref{maint} to the \emph{singular} case $p<2$.\\

\begin{figure}[ht]
 \centering
\includegraphics[trim={3cm 8cm 2cm 7cm},clip, scale=0.8]{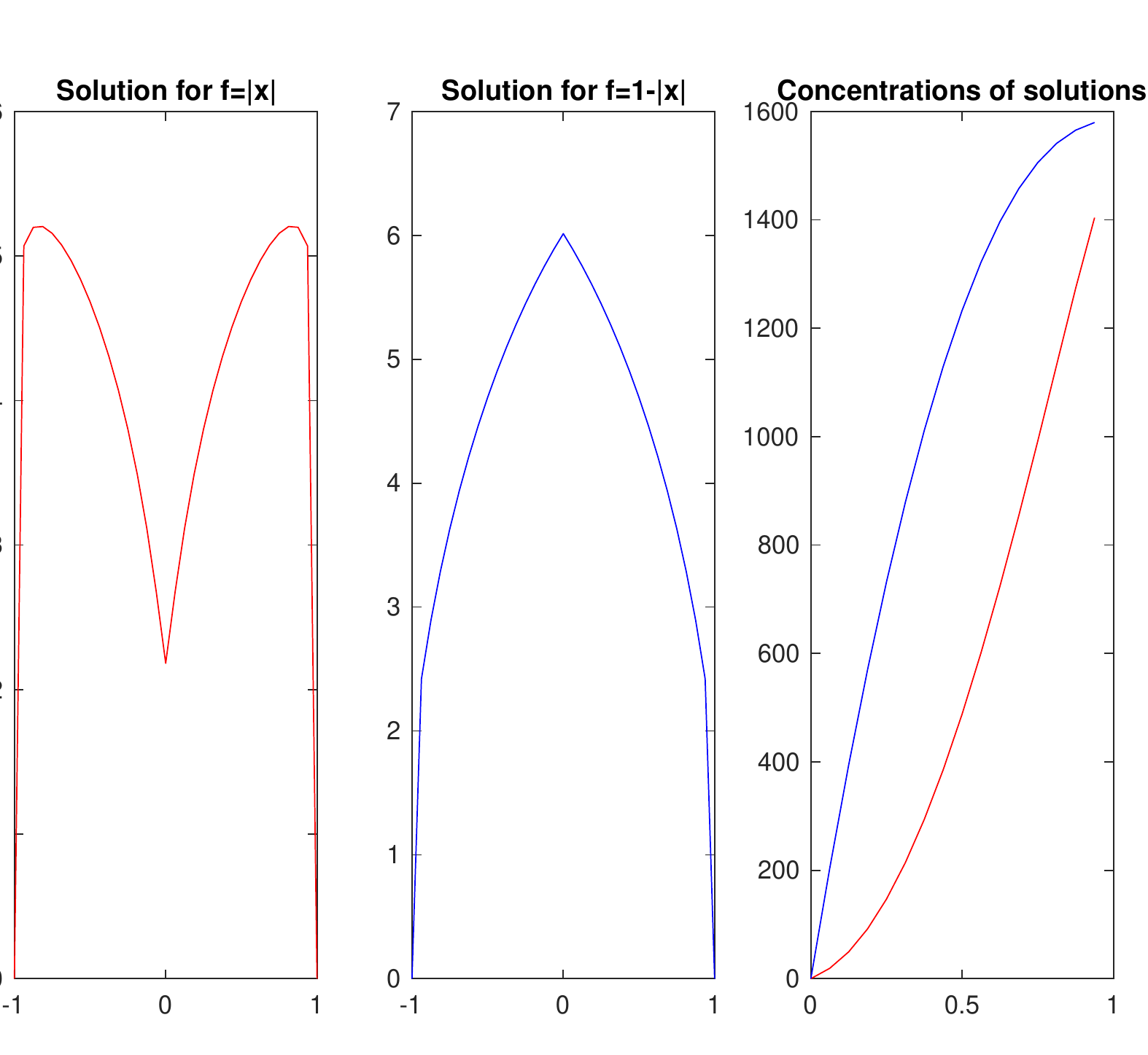}
\caption{From left to right: plot of $u$ with the choice $f=|x|$, plot of $v$ and comparison of mass concentrations of $p-1$ powers of $u^{\#}$ and $v$.}
 \label{fig:useCase}
\end{figure}

\noindent {$\bullet$} A natural question would be to ask comparing the solution $u$ to problem \eqref{mainpp} with the radial solution $v$ to the \emph{nonlinear} problem
\begin{equation} \label{maip}
\left\{
\begin{array}[c]{lll}%
(-\Delta_{p})^{s}v=f^{\#} & & \text{in }%
\Omega^{\#},\\
\\
v=0 & & \text{on }\partial\Omega^{\#}.
\end{array}\right.
\end{equation}
Actually from inequality \eqref{mainineqnon}, the radiality of $v$ would lead quite easily to the integral comparison
\begin{equation}\label{openpr}
\int_{B_r}\int_{B_r^c}\frac{|u^{\#}(x)-u^{\#}(y)|^{p-1}}{|x-y|^{N+sp}}\,dx\,dy
 \leq
\int_{B_r}\int_{B_r^c}\frac{|v(x)-v(y)|^{p-1}}{|x-y|^{N+sp}}\,dx\,dy.
\end{equation}
An interesting open problem would be to apply maximum principle arguments to inequality \eqref{openpr} in order to derive an integral comparison between the $p-1$ powers of $u,\,v$ in the sense of \cite{fermess}, namely an inequality of the type
\begin{equation}
\int_{B_{r}(0)}(u^{\#})^{p-1}dx\leq
\int_{B_{r}(0)} v^{p-1}dx,\quad r>0.\label{compconcp}
\end{equation}
Though such result would appear quite natural, at the same time it seems very difficult to prove: indeed, it seems to require some arguments that are fairly different to the ones established in \cite{FeroneVolzone}. Nevertheless, the following numerical simulation suggests that \eqref{compconcp} is quite natural to expect. The plots in Figure \ref{fig:useCase}, obtained by a suitable implementation of the robust numerical methods in \cite{TesoLind} seem to confirm our guess. Indeed, we have considered problem \eqref{mainpp} in the case $N=1$, $p=3$, $s=\frac12$, $\Omega=(-1,1)$, and we have compared, in terms of mass concentration, the solution $u$ when $f(x)=|x|$ with the solution $v$ to problem \eqref{maip} when $f^\#(x)=1-|x|$. The theoretical question will be an object of future investigations.

\newpage
\section*{Acknowledgments}

V.F. was partially supported by Italian MIUR through research project PRIN 2017 ``Direct and inverse problems for partial differential equations: theoretical aspects and applications''. B.V. was partially supported by Gruppo Nazionale per l'Analisi Matematica, la Probabilit\`a e le loro Applicazioni (GNAMPA) of Istituto Nazionale di Alta Matematica (INdAM). Both authors are members of GNAMPA of INdAM. B.V. wishes to warmly thank L. Brasco for fruitful discussions and valuable suggestions.
\nc


%
2000 \textit{Mathematics Subject Classification.}
35B45,  
35R11,   	
35J25. 


%
\textit{Keywords and phrases.} Symmetrization, fractional Laplacian,
 nonlocal elliptic equations.

\end{document}